\documentclass[11pt]{article}
\usepackage{graphicx}     
\usepackage{natbib}       
\usepackage{enumerate}
\usepackage{cite}
\usepackage{subfigure}
\usepackage{epstopdf}
\usepackage[utf8]{inputenc}
\usepackage[T1]{fontenc}
\usepackage{amsfonts,latexsym,amssymb,amsmath}
\begin{document}

\title{Control of Pitch-Flap Instabilities in Helicopter Rotors using Delayed Feedback} 
\author{Rudy Cepeda-Gomez\\\emph{a+i Engineering, Berlin, Germany}\\ \texttt{rudy.cepeda-gomez@a-{}-i.com}}
\date{}
\maketitle

\begin{abstract}
The problem of vibration suppression in helicopter rotors is the main topic of this paper. We revisit a strategy for the stabilization of such systems based on delayed feedback. By means of an improved analysis, made possible by recently developed mathematical tools, we present insights on how the delay can be used as stabilization tool for this problem. The results are supported with numerical simulations. For simplicity, only the case of hovering flight is considered.
\end{abstract}

\section{Introduction}
The main rotor is the principal source of vibrations in a helicopter and the reduction of these vibrations is important to increase the life of the components of the airframe as well as the comfort of passengers and crew. A considerable research effort has been devoted to the study of active control systems for helicopter rotors, in order to guarantee the stability of its motion. Most works on this topic fall in one of two categories: Higher Harmonic Control (HHC) or Individual Blade Control (IBC). Works in the HHC category, like those of \citet{Hall1995,Mura2014,Mura2014-2} consider the control forces to be applied to the swashplate by means of actuators. The IBC works \citep{Bittanti2002,Arcara2000,King2014,Shen2006} consider that each blade can be actuated independently. The main difference between the two approaches is the number of available degrees of freedom for the control action. Many control algorithms, nevertheless, can be implemented with either actuation technique. \citet{Friedmann1995} compare several different control methods in both categories. A more recent survey of the developments in these topics was presented by \citet{Friedmann2014}.

The present paper revisits a control strategy presented by \citet{krodkiewski2000}, which can be implemented in both HHC and IBC frameworks. This control logic uses delayed feedback to improve the stability of the motion of the helicopter rotor. Delayed feedback has been previously used for the stabilization of periodic motion by \citet{Pyragas1992} and for the suppression of vibrations by \citet{Olgac1994}. Using novel analysis techniques, we challenge the published results and show a corrected stability region in the parametric space.

Although the original works of \citet{Faragher1996} and \citet{krodkiewski2000} consider both the cases of hovering and horizontal flight, we restrict the scope of this work to the case of vertical motion only. In this way, we can apply tools for the stability analysis of Linear Time Invariant (LTI) systems affected by delays. The case of horizontal flight, and the study of the Linear Time Periodic (LTP) systems with delays, is a matter of further research.

The paper is organized as follows. Section \ref{sec:dyn} introduces the dynamic model of the rotor in hovering flight, borrowed from \citet{Bramwells}. Section \ref{sec:unctrl} presents the instabilities that may appear in the operation of an uncontrolled rotor. Section \ref{sec:crtl} introduces the control law and the stability analysis. Some concluding remarks are finally presented in section \ref{sec:conc}.

\section{Dynamics of the Rotor}\label{sec:dyn}
\begin{figure}[!b]
\centering
\includegraphics[scale=0.5]{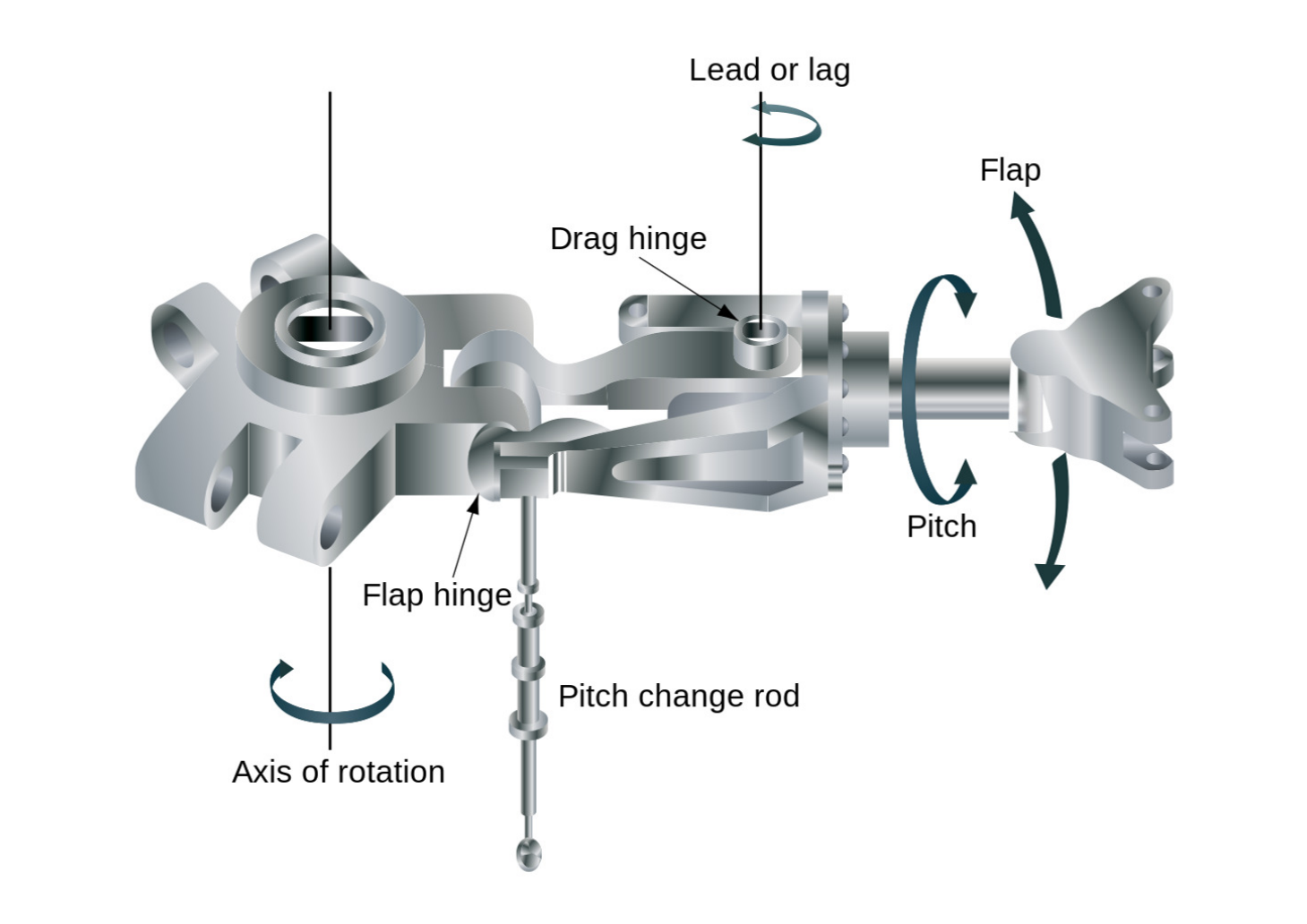}
\caption{Typical hinge arrangement and degrees of freedom of a helicopter rotor blade. (Source: Wikimedia Commons)}
\label{fig:rotorhead}
\end{figure}

A helicopter rotor blade is mounted on a set of hinges which allow three angular degrees of freedom. A typical arrangement of the hinges is shown in figure \ref{fig:rotorhead}. The \emph{flapping} motion is defined as an up and down rotation in a plane which contains both the blade and the shaft. The flapping angle, $\beta$, is defined as positive when the blade moves upwards. A flapping blade rotating at high speeds is subject to large Coriolis moments in the plane of rotation, and the \emph{lag} hinge is introduced to alleviate these moments. This hinge allows the motion of the blade in the same plane of rotation. The lagging angle, $\xi$, is positive when the blade moves in the same sense of the shaft rotation. The \emph{pitching} or \emph{feathering} motion, denoted by $\theta$, is a rotation about an axis parallel to the blade span.

The pitch-flap flutter is an unstable behavior produced by the coupling of the pitching and flapping motions of the rotor. It can lead to very high oscillatory loads in the pitch control mechanism. The following dynamic model, adapted by \citet{krodkiewski2000} from the work of \citet{Bramwells}, describes the coupling between $\theta$ and $\beta$ for an uncontrolled rotor in hovering flight
\begin{equation}
M\left[\begin{array}{c}\ddot{\theta}(\psi)\\\ddot{\beta}(\psi)\end{array}\right]+C\left[\begin{array}{c}\dot{\theta}(\psi)\\\dot{\beta}(\psi)\end{array}\right]+K\left[\begin{array}{c}\theta(\psi)\\\beta(\psi)\end{array}\right]=0,
\label{eq:fluttermodel}
\end{equation}
with the inertia, damping, and stiffness matrices:
\begin{equation}
M=\left[\begin{array}{cc}1&-12\frac{r_g\sigma}{c_h}\\0&1\end{array}\right],\ 
C=\left[\begin{array}{cc}\frac{\gamma}{8}&0\\ 0&\frac{\gamma}{8}\end{array}\right],\ 
K=\left[\begin{array}{cc}\nu_1^2&-12\frac{r_g\sigma}{c_h}\\-\frac{\gamma}{8}&\lambda _1^2\end{array}\right]
\label{eq:matrices}
\end{equation}
As it is customary in helicopter dynamics, the independent variable in \eqref{eq:fluttermodel} is not the time $t$, but the azimuth angle $\psi=\Omega t$. Here $\Omega$ is the angular speed of the rotor shaft, which is assumed constant during the operation of the aircraft. In \eqref{eq:fluttermodel} a dot over a variable denotes its derivative with respect to $\psi$.

It is important to state that in the derivation of \eqref{eq:fluttermodel} it is assumed that no lag motion occurs.

Among the parameters presented in \eqref{eq:matrices}, the stability is usually studied with respect to the two parameters which have the greatest influence on it. These parameters are $\sigma$, which is the distance of the center of gravity (c.g.) of the blade from the center of pressure, expressed as a fraction of the chord; and the non rotating torsional natural frequency $\nu_1$. The other parameters, described in table \ref{tb:params}, are taken as fixed. Notice that the natural frequencies $\nu_1$ and $\lambda_1$ are given as fractions of $\Omega$. Table \ref{tb:params} also shows the values used for the fixed parameters in the numerical analysis section. These values are taken from the work of \citet{Faragher1996}.
\begin{table}[!b]
\centering
\caption{Parameters of the Dynamic Model}
\label{tb:params}
\begin{tabular}{clc}
\hline \hline
Param.&Description&Value\\
\hline \hline
$r_g$&Span wise location of blade's cg&4.1 m\\
$c_h$&Chord&0.527 m\\
$\gamma$&Lock Number&6.95\\
$\lambda_1$&First flap natural frequency&1.1\\
$\sigma$&Pos. of the cg from the center of pressure&$[0,0.08]$\\
$\nu_1$&Torsional natural frequency&$[0,\sqrt{8}]$\\ 
\hline \hline
\end{tabular}
\end{table}

Using two sticks in the cockpit the pilot can exert a control moment which directly affects the pitch motion of the rotor. The basic effect of a pitch change is a change of the lift force produced by the blade. If all blades change the pitch simultaneously, the lift of the whole rotor is changed and vertical motion ensues. This is the so called collective pitch control and is changed using one stick in the cockpit. If the pitch of the blades changes as a function of its angular position, a directional force can be created and horizontal flight (forwards, rearwards, sideways) is attained. This is the cyclic pitch control and is controlled by a second stick in the cockpit. The uncontrolled model \eqref{eq:fluttermodel} is then expanded using an input signal $u(\psi)$, which corresponds to the motion of the lower end of an actuation rod. Following again \citet{Faragher1996}, the controlled  model with the parametric values under consideration becomes  
\begin{equation}
M\left[\begin{array}{c}\ddot{\theta}(\psi)\\\ddot{\beta}(\psi)\end{array}\right]+C\left[\begin{array}{c}\dot{\theta}(\psi)\\\dot{\beta}(\psi)\end{array}\right]+K\left[\begin{array}{c}\theta(\psi)\\\beta(\psi)\end{array}\right]=Bu(\psi)=\left[\begin{array}{c}3777\nu_1^2\\0\end{array}\right]u(\psi).
\label{eq:ctrlmodel}
\end{equation}

Equations \eqref{eq:fluttermodel} and \eqref{eq:ctrlmodel} are models for hovering or vertical flight only. In this situation, the matrices in \eqref{eq:matrices} are constant, and we deal with linear, time-invariant (LTI) systems. When forward (or sideways) motion of the helicopter is considered, the aerodynamic forces and moments become time variant. More specifically, in these cases the matrices in \eqref{eq:matrices} are periodic, with a time period $T=2\pi\Omega$ or simply $T=2\pi$ when the azimuth angle $\psi$ is used as independent variable. This case is certainly more complicated and is left aside for future studies.

\section{Stability Analysis of the Uncontrolled System}\label{sec:unctrl}
With the state vector $x\left(\psi\right)=\left[\theta\,\beta\,\dot{\theta}\,\dot{\beta}\right]^T$, equation \eqref{eq:fluttermodel} can be represented in state space as the fourth order system
\begin{equation}
\dot{x}(\psi)=A_ux(\psi)=\left[\begin{array}{cc}0_{22}&I_{22}\\-M^{-1}K&-M^{-1}C\end{array}\right]x(\psi),
\label{eq:ssuctrl}
\end{equation}
with $0_{22}$ being a 2-by-2 matrix with all elements equal to zero and $I_{22}$ the identity matrix of the same dimension.

As the parameters $\sigma$ and $\nu_1^2$ change while keeping the other constant, two types of instability can appear. When $A_u$ has a purely real unstable root, the unstable mode is known as \emph{pitch divergence} and it is evidenced by a steady increase in the pitch and flap angles. When the unstable roots are complex conjugate, the oscillatory instability known as \emph{pitch-flap flutter} appears.

The boundary at which the transition from stability to pitch divergence occurs is given by the presence of a characteristic root at the origin of the complex plane. A sufficient condition for $A_u$ to have a root at the origin is that the determinant of the stiffness matrix is equal to zero, i.e., $|K|=0$. From here we obtain
\begin{equation}
\nu_1^2=\frac{3\gamma r_g\sigma}{2c_h\lambda _1^2}
\label{eq:divbound}
\end{equation}

To establish the boundary between stability and flutter, it is necessary to find a conjugated pair of purely imaginary characteristic roots at a frequency $\omega_f$. This is achieved by solving the characteristic equation $|j\omega_f I_{44}-A_u|=0$, leading to the following two equations in terms of the crossing frequency
\begin{subequations}
\label{eq:flutbound}
\begin{align}
\nu_1^2&=2\omega_f^2-\lambda_1^2,\\
\sigma&=\frac{c_h\left(\gamma^2\omega_f^2+64\omega_f^4-128\omega_f^2\lambda_1^2+64\lambda_1^4\right)}{96r_g\gamma\left(\omega_f^2-1\right)}.
\end{align} 
\end{subequations}

\begin{figure}
\centering
\includegraphics[scale=0.8]{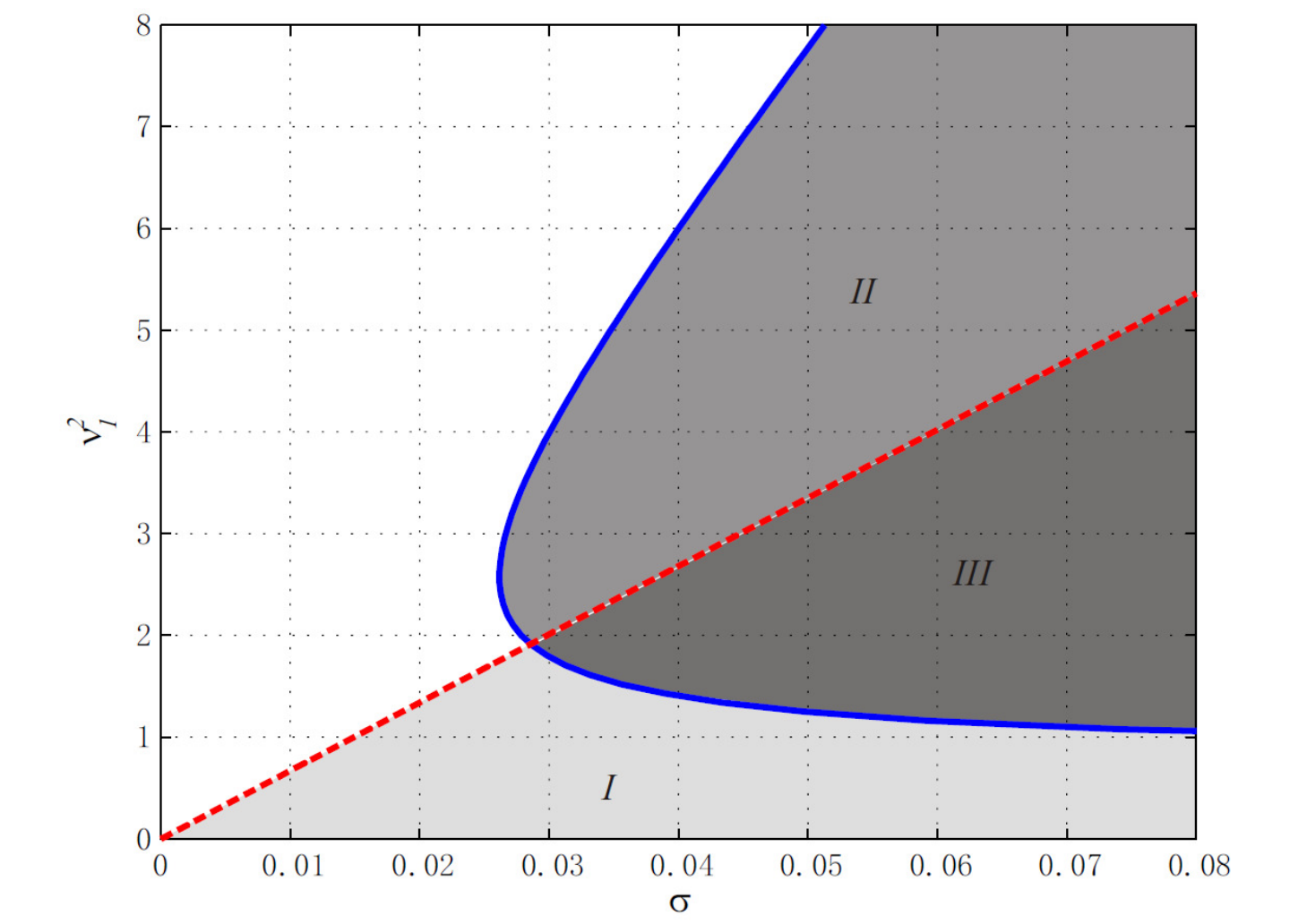}
\caption{Pitch-flap stability boundaries for an uncontrolled system. The non shaded region represents stable parametric combinations. The shaded regions indicate pitch divergence (I), pitch-flap flutter (II) and both (III) unstable behaviors.}
\label{fig:unctrl}
\end{figure}

The two stability boundaries are presented in Fig.~\ref{fig:unctrl}. The red, dashed line depicts the boundary at which an unstable pole introduces pitch divergence, whereas the blue, solid line indicates a change in stability to the pitch-flap flutter. These two lines divide the parametric space in four zones. The region without shading corresponds to stable parametric selections. The light shaded zone at the bottom of the plot (zone I) indicates the presence of a pitch divergence mode. The darker shaded region at the top right (zone II) indicates pitch-flap flutter behavior, and the darkest region to the right (zone III) corresponds to the worst possible conditions, when both unstable modes are present.
\section{Stability Analysis of the Controlled System}\label{sec:crtl}
\citet{Faragher1996} and \citet{krodkiewski2000} present a method to increase the size of the stable region in Fig.~\ref{fig:unctrl} by means of delayed feedback. Based on the work of \citet{Pyragas1992}, the authors of these works propose the control law
\begin{equation}
u\left(\psi\right)=a\left(\beta\left(\psi\right)-\beta\left(\psi-\tau\right)\right)+b\left(\dot{\beta}\left(\psi\right)-\dot{\beta}\left(\psi-\tau\right)\right),
\label{eq:control}
\end{equation}
with $\tau=2\pi$. The idea behind \eqref{eq:control} is to use as error signal the difference between the current state and its value one period earlier. \citet{Faragher1996} and \citet{krodkiewski2000} claim that making the delay exactly equal to the period of the motion greatly increases the region of asymptotic stability of the uncontrolled system. On the contrary, this particular selection of the delay value actually introduces instability into the system. In the following paragraphs, we analyze the effect of the delay on the stability of the system and show under which conditions the delayed feedback is beneficial or not.

Under \eqref{eq:control} the controlled system \eqref{eq:ctrlmodel} is expressed as
\begin{equation}
\dot{x}\left(\psi\right)=Ax\left(\psi\right)+A_dx\left(\psi-\tau\right)
\label{eq:delay}
\end{equation}
with
\begin{align*}
A&=\left[\begin{array}{cc}0_{22}&I_{22}\\-M^{-1}K+E&-M^{-1}C+F\end{array}\right],\quad A_d=\left[\begin{array}{cc}0_{22}&0_{22}\\-E&-F\end{array}\right],\\
E&=\left[\begin{array}{cc}0&3777v_1^2a\\0&0\end{array}\right],\quad F=\left[\begin{array}{cc}0&3777v_1^2b\\0&0\end{array}\right].
\end{align*}

The stability analysis of a Linear Time Invariant, Time Delay System (LTI-TDS) like \eqref{eq:delay} can be done following the two propositions of the Cluster Treatment of Characteristic Roots presented by \citet{Olgac2002}. These two propositions state the following: 
\begin{itemize}
\item Although an LTI-TDS has infinitely many characteristic roots, which change as the delay does, there is only a finite number of possible purely imaginary characteristic roots, which appear at periodically spaced values of the delay. That is, if $s=j\omega_c$ is a characteristic root of the system for $\tau=\tau_c$, the same imaginary root appears when $\tau=\tau_c+2k\pi/\omega_c$, $k=1,\,2,\,3\,\ldots$.
\item When the delay increases an imaginary characteristic root can move to the left or to the right of the complex plane, stabilizing or destabilizing the system. This transition is invariant for roots created by the periodicity. That is, if by increasing the delay from $\tau_c$ to $\tau_c+\varepsilon$, $0<\varepsilon<<1$, the root at $s=j\omega_c$ moves to the left, it does the same when the delay increases from $\tau=\tau_c+2\pi/\omega_c$ to $\tau=\tau_c+2\pi/\omega_c+\varepsilon$.
\end{itemize}

According to these ideas, to determine the values of the delay for which the system is stable we must find the purely imaginary characteristic roots of system \eqref{eq:delay}. For this, we form the characteristic equation
\begin{equation}
\det\left(sI_{44}-A-A_de^{-\tau\,s}\right)=0.
\label{eq:ce}
\end{equation}
Since $A_d$ has rank 1, this equation has the form
\begin{equation}
P(s)+Q(s)e^{-\tau\,s}=0.
\label{eq:cepq}
\end{equation}
By replacing $s=j\omega_c$ in \eqref{eq:cepq}, and reorganizing we obtain
\begin{equation}
\frac{P\left(\omega_c\right)}{Q\left(\omega_c\right)}=-e^{-j\omega_c\,\tau}.
\label{eq:cepq2}
\end{equation}
After equating the magnitudes of both sides of \eqref{eq:cepq2} we arrive to an expression in terms of $\omega_c$
\begin{equation}
\left|P\left(\omega_c\right)\right|^2-\left|Q\left(\omega_c\right)\right|^2=0.
\label{eq:mag}
\end{equation}
For the particular case at hand, equation \eqref{eq:mag} is an eighth degree polynomial in $\omega_c$, with nonzero coefficients only for the even exponents of $\omega_c$. Introducing a new variable $\eta=\omega_c^2$, we get an expression of the form
\begin{equation}
\eta^4+c_3\eta^3+c_2\eta^2+c_1\eta+c_0=0.
\label{eq:omega2}
\end{equation}
The solutions of \eqref{eq:omega2} with $\eta\geq0$ correspond to valid crossing frequencies $\omega_c=\pm\sqrt{\eta}$ for the LTI-TDS system \eqref{eq:delay}. This equation shows that there can be up to 4 crossing points at the imaginary axis. To find the delay at which such crossing frequencies occur, we equate the phase angles of both members of \eqref{eq:cepq2} and solve for $\tau_c$
\begin{equation}
\tau=\frac{1}{\omega_c}\left(\angle Q\left(\omega_c\right)-\angle P\left(\omega_c\right)+\pi\right).
\label{eq:tau}
\end{equation}

The explicit versions of \eqref{eq:omega2} and \eqref{eq:tau} are very long, and thus left outside of this document.

To find whether an imaginary root crosses to the left or to the right of the complex plane, the root tendency for each $\omega_c$ and one of its corresponding $\tau$ is found using its definition \citep{Olgac2002}
\begin{equation}
RT=\text{signum}\left(\Re\left[\frac{ds}{d\tau}\right]\right).
\label{eq:rt}
\end{equation}
The root sensitivity $ds/d\tau$ is found from \eqref{eq:ce} using implicit differentiation. When $RT=1$, the root crosses to the right half of the complex plane, which introduces instability. If $RT=-1$ the root moves to the stable left half of the complex plane.

With the help of \eqref{eq:omega2}, \eqref{eq:tau} and \eqref{eq:rt}, we can determine for which values of the delay the system is stable given a certain parametric selection. The procedure is simple. After finding the crossing frequencies and the delays that generate them we create a table of stability transition points. Starting with the number of unstable characteristic roots of the non delayed system (which is equivalent to the number of eigenvalues of $A$ with positive real part), we increase $\tau$ until a crossing delay is reached. If $RT=1$ for that point, we increase the number of unstable roots by 2, or decrease it if $RT=-1$. Those intervals of $\tau$ for which the number of unstable roots is zero are the delays that guarantee stable operation.

\subsection{Stabilization of Pitch Divergence}
\citet{Faragher1996} states that the control action \eqref{eq:control} cannot stabilize a system that presents pitch divergence. They arrived to this conclusion by means of a parametric search. Here we provide an analytic proof of it.

When the helicopter rotor presents pitch divergence, the non delayed case has a real positive characteristic root. For this root to become stable, a stabilizing transition at the origin of the complex plane must take place. That is, we need a crossing frequency $\omega=0$ for which $RT=-1$. For \eqref{eq:omega2} to have a root at $\omega=0$, the independent term $c_0$ must be zero. Considering the explicit form of this coefficient, the condition for this to happen is
\begin{equation}
c_0=\lambda_1^2\left(\lambda_1^2-944\gamma a\right)\nu_1^4+\frac{\gamma r_g\sigma}{c_h}\left(3\lambda_1^2-1416\gamma a\right)\nu_1^2+\frac{9}{4}\frac{\gamma^2r_g^2\sigma^2}{c_h^2}=0.
\label{eq:c0}
\end{equation}
This equation shows that for $\sigma,\,\nu_1^2>0$ it is impossible to set the control parameters to a value such that $c_0=0$. Therefore no crossing of the imaginary axis occurs at the origin (in either direction) and the transition from pitch divergence to stability cannot take place.
\subsection{Stabilization of Pitch-Flap Flutter}
In order to show that the proposed control law \eqref{eq:control} is able to eliminate the pitch flap instability in a wide range of values of the operating parameters, \citet{Faragher1996} explores the complex plane\footnote{We use $\chi=\Re(s)$ to avoid confusion with the parameter $\sigma$.} $s=\chi+j\omega$ searching for locations where the loci of zero points of the real and imaginary parts of the characteristic equation \eqref{eq:ce} intersect. Their variation of parameters follows this logic: after selecting the values of $\sigma$ and $\nu_1$, the control gains $a$ and $b$ are varied, searching for the optimal location of the characteristic roots, i.e., the location where the rightmost roots have the smallest possible real part. 

\begin{figure}
\centering
\subfigure[$$ Original search region.]{\includegraphics[scale=0.3]{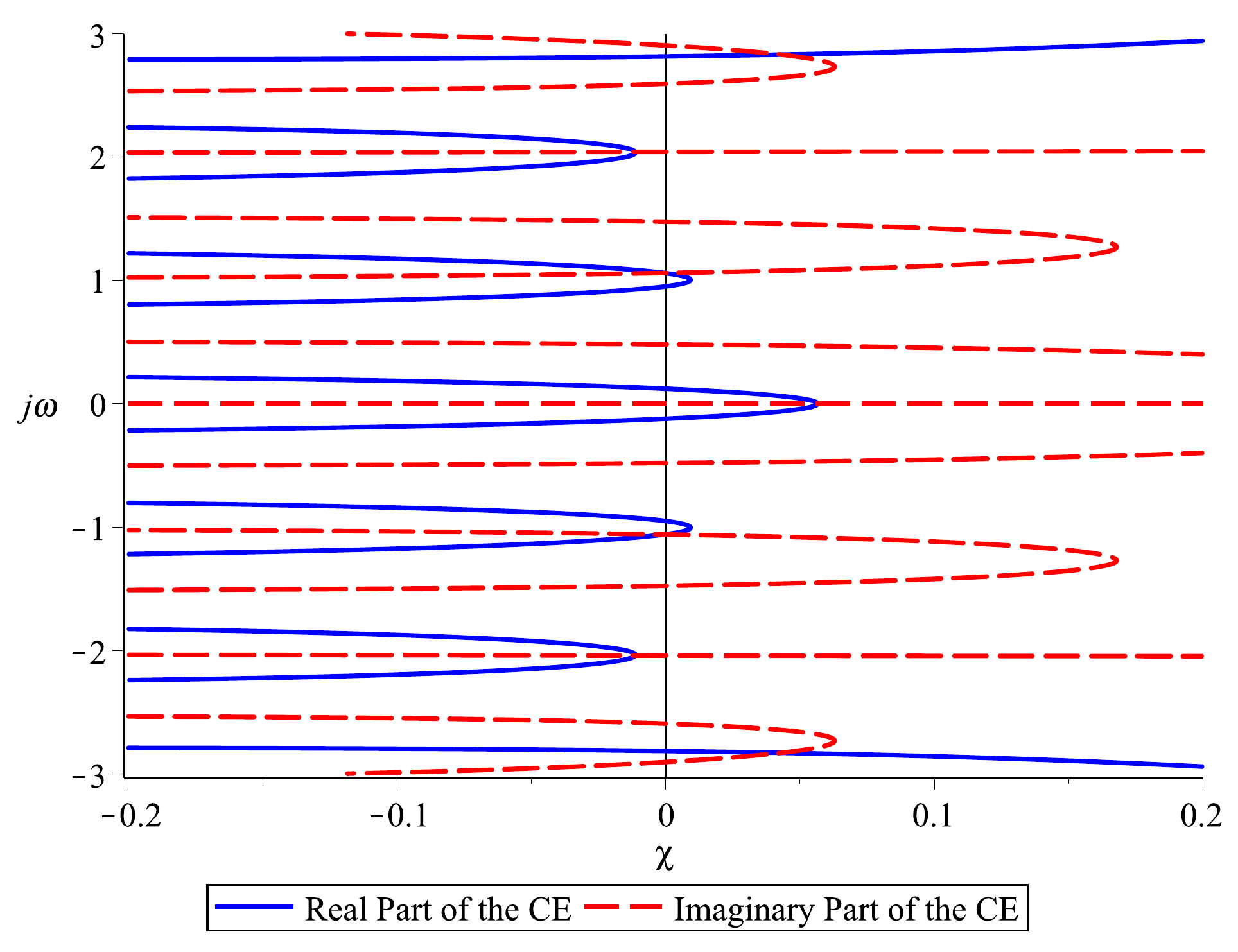}\label{fig:rlocus1}}
\subfigure[$$ Expanded search region.]{\includegraphics[scale=0.3]{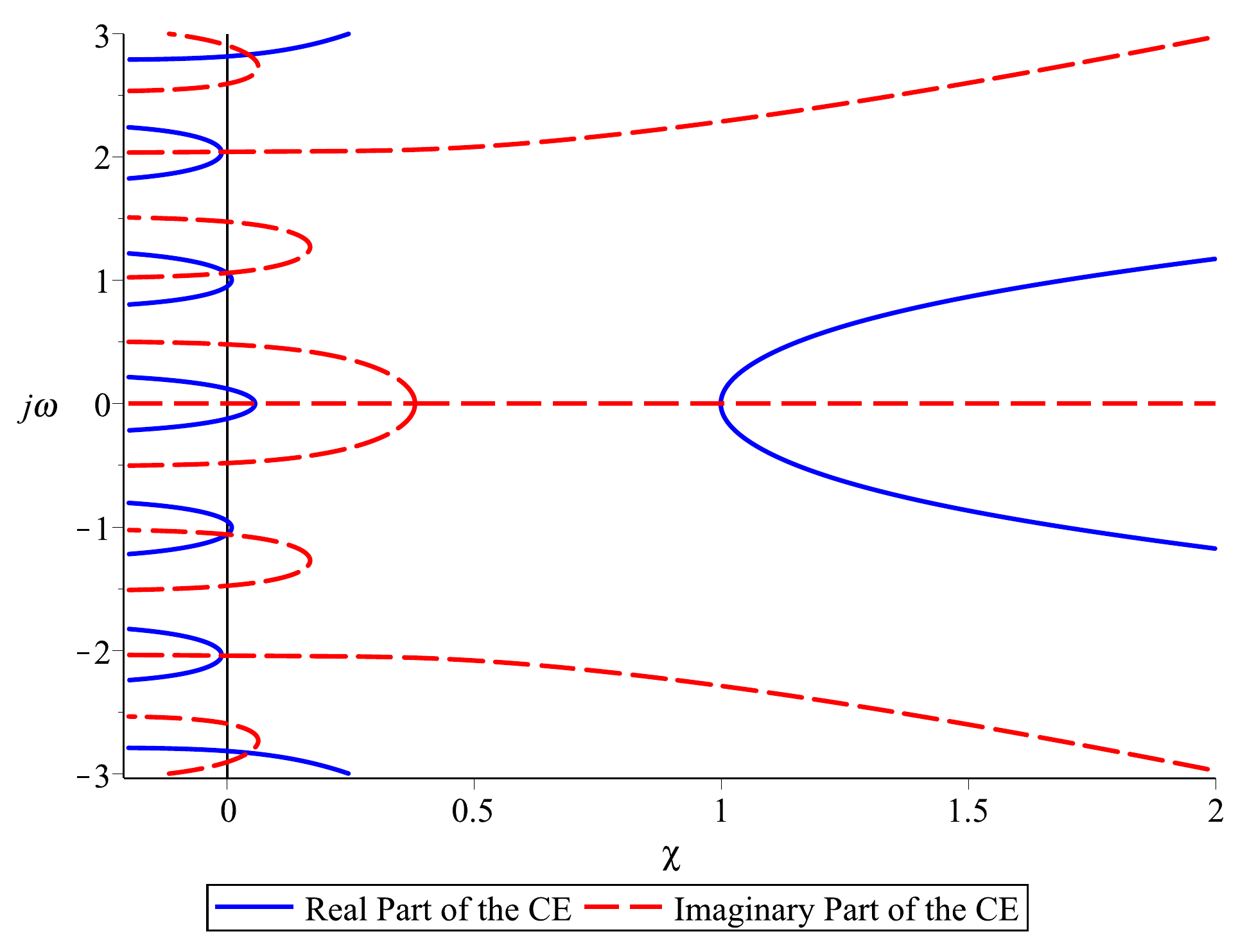}\label{fig:rlocus2}}
\caption{Root locus of zero points of the real (blue, solid) and imaginary (red, dashed) parts of the characteristic equation \eqref{eq:ce} for $\sigma=0.08$, $\nu_1^2=10.8$, $a=6.75\times 10^{-4}$, $b=0.6\times 10^{-4}$, $\tau=2\pi$. The intersections indicate the location of characteristic roots of the system.}
\end{figure}

This method has a big drawback: since \eqref{eq:ce} is a transcendent equation it has infinitely many roots and any parametric exploration would never be able to find all of them to determine which one is the rightmost. We used the same searching method they applied and did not obtain the same root loci and characteristic roots presented by \citet{krodkiewski2000} in figure 5 and by \citet{Faragher1996} in figure 6.19. For the parametric selections stated in those works\footnote{Please note that in this paper we are not using the $10^{-4}$ scaling factor in the control gains introduced by \citet{Faragher1996} or \citet{krodkiewski2000}.} and considering the same search region, we obtained the plot in Fig.~\ref{fig:rlocus1}. The presence of three unstable characteristic roots is clear. Furthermore, one of these roots is real, implying that the control action is introducing pitch divergence into a parametric selection that did not exhibit this behavior. The roots in Fig.~\ref{fig:rlocus1} are, however, not rightmost. If the search region is increased further to the right, another real root shows up, as it is shown in Fig.~\ref{fig:rlocus2}. This highlights how parametric searches of this kind may miss very important information.

\begin{figure}[tb]
\centering
\subfigure[$$ $\tau=2\pi$]{\includegraphics[scale=0.55]{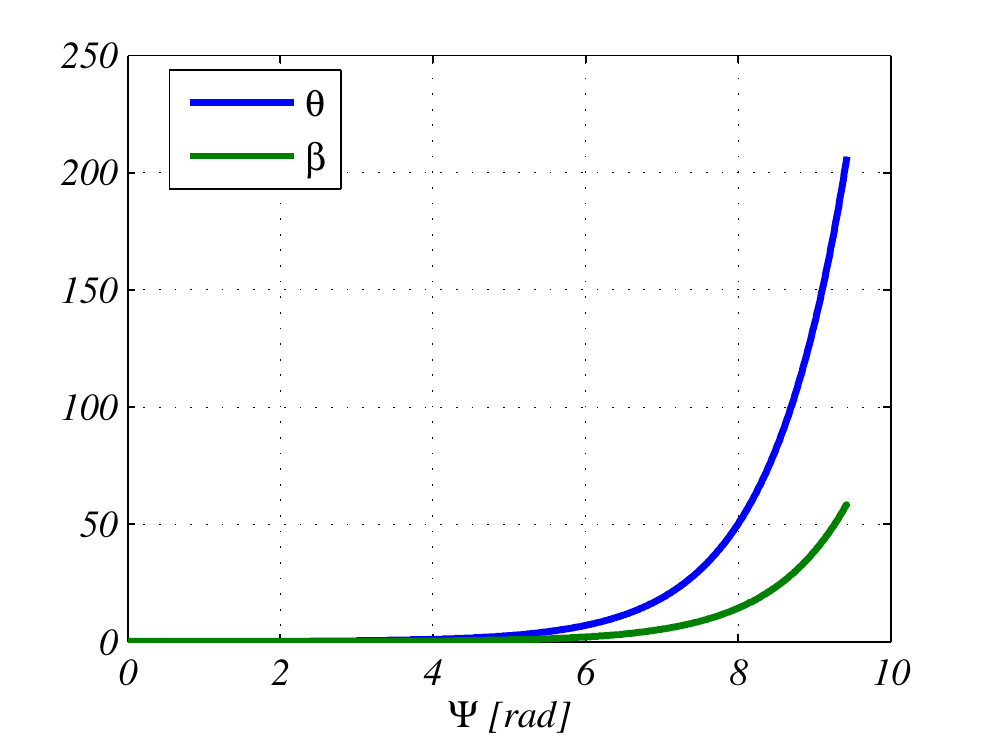}\label{fig:sim1}}
\subfigure[$$ $\tau=0.2296$]{\includegraphics[scale=0.55]{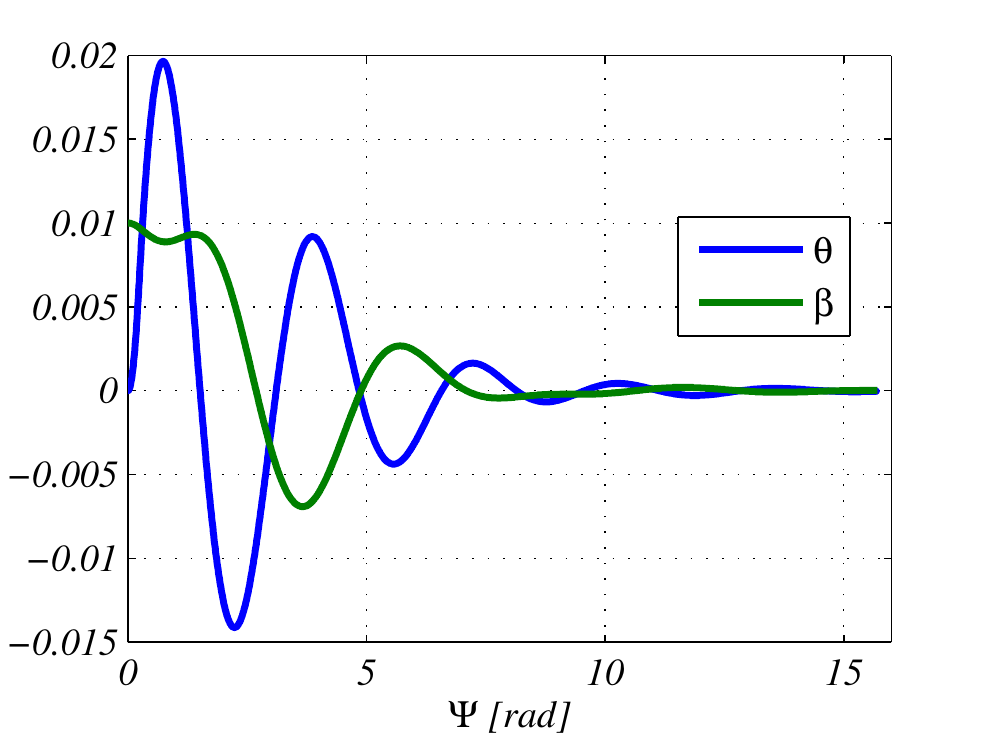}\label{fig:sim2}}
\caption{Simulation results for a system with $\sigma=0.08$, $\nu_1^2=10.8$, $a=6.75\times 10^{-4}$, $b=0.6\times 10^{-4}$ with two different delays.}
\end{figure}

In order to verify the results of Fig.~\ref{fig:rlocus2}, we used the QPmR function \citep{vyhlidal2009} to numerically find the roots in the same region of the complex plane. The results coincide: QPmR finds the roots at  $0.9927$, $0.0595$, $-8\times 10^{-4}+j1.0623$, $-1.15\times 10^{-2}+j2.0439$, $4.10\times 10^{-2}+j2.8268$ and $-0.196+j3.5967$. We also simulated the system using the same initial condition stated by \citet{krodkiewski2000}: $[\theta\,\beta]^T=[0\,0.01]^T$. The simulation results of Fig.~\ref{fig:sim1} show a time constant which is consistent with the presence of a real unstable root located near to $s=1$. 

For the selected control parameters in the previous example, which are declared as the optimal, the delay value selected by \citet{krodkiewski2000} actually creates instability. By applying the CTCR methodology described at the beginning of the section, we can find which values of the delay are actually able to stabilize the system. For this particular parametric selection, the solution of \eqref{eq:omega2} gives three crossing frequencies, which are listed in Table \ref{tb:crossings1} together with their corresponding delays and root tendencies.
\begin{table}
\centering
\caption{Crossing Frequencies and Delays for $\sigma=0.08$, $\nu_1^2=10.8$, $a=6.75\times 10^{-4}$, $b=0.6\times 10^{-4}$}
\label{tb:crossings1}
\begin{tabular}{ccc}
\hline \hline
$\omega_{ci}$&$RT$&$\tau_{ci0}$\\
\hline \hline
2.1949&$-1$&0.0852\\
1.0525&1&0.3580\\
3.0268&1&1.519\\
\hline \hline
\end{tabular}
\end{table} 

\begin{figure}[tb]
\centering
\includegraphics[scale=0.4]{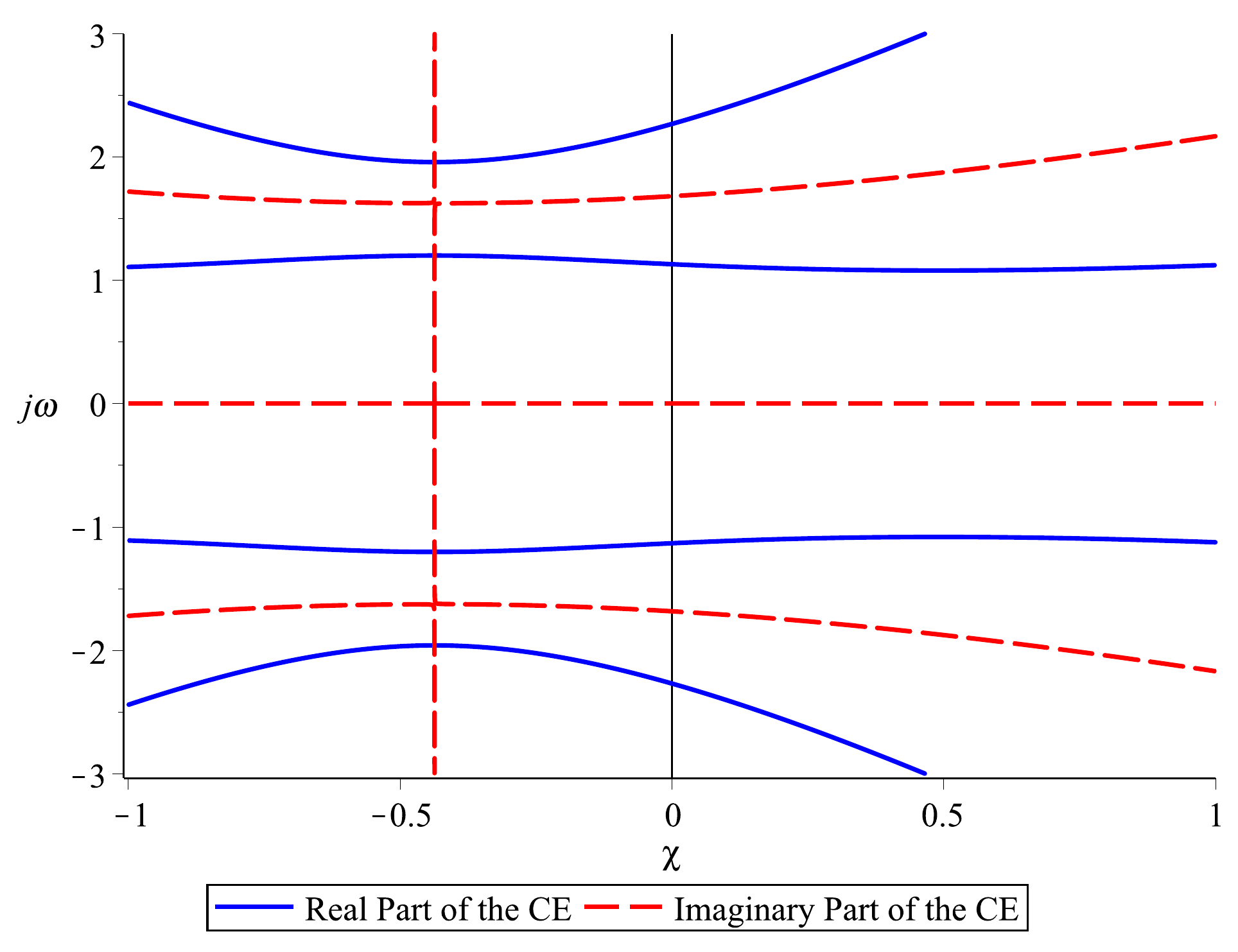}
\caption{Root locus of the real (blue, solid) and imaginary (red, dashed) parts of the characteristic equation \eqref{eq:ce} for $\sigma=0.08$, $\nu_1^2=10.8$, $a=6.75\times 10^{-4}$, $b=0.6\times 10^{-4}$, $\tau=0.2296$. The intersections indicate the location of characteristic roots of the system.}
\label{fig:optimum}
\end{figure}

The smallest crossing delay has $RT=-1$. This implies that the system, which presents pitch-flap flutter instability for $\tau=0$, recovers its stability when the delay is larger than $\tau=0.0852$. However, it loses its stability for delays larger than $\tau=0.3560$. The optimal delay value is then in the interval $[0.0852\ 0.3560]$. By varying the delay whilst keeping the control gains fixed, which is a numerical search approach similar to that used by \citet{krodkiewski2000}, we found the optimal delay to be $\tau=0.2296$, as it is shown in Fig.~\ref{fig:optimum}. The rightmost characteristic roots for this case are $-0.4368+j1.2018$ and $-0.4368+j1.9596$. Figure \ref{fig:sim2} presents the simulation results for this case. The search of the optimal delay value is made possible by the knowledge of the exact stability interval with respect to the delay. Knowing the upper and lower bounds, we can search in a finite space where we are sure to find the rightmost root with smallest possible real part.

\begin{figure}[tb]
\centering
\includegraphics[scale=0.7]{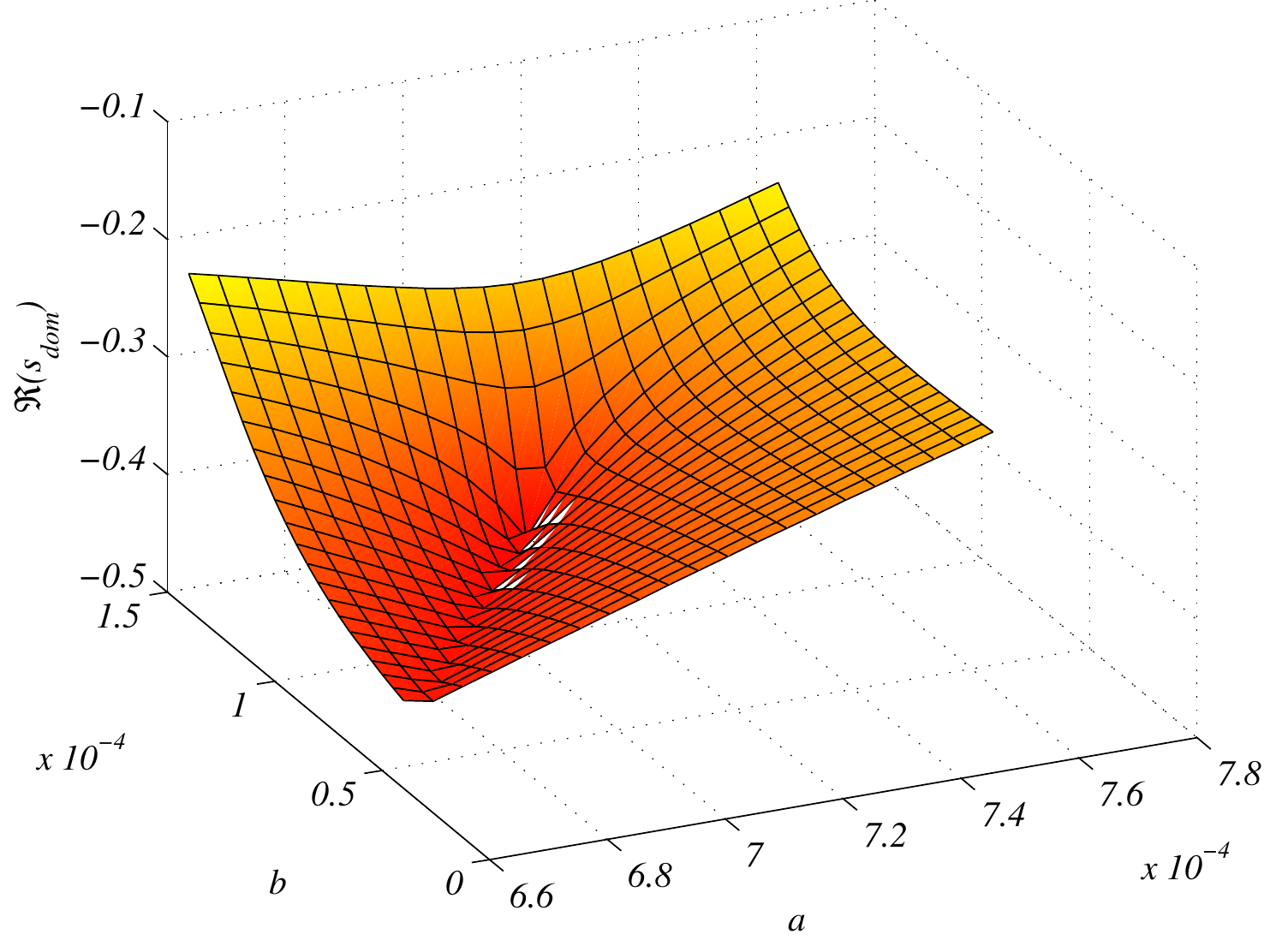}
\caption{Real part of the rightmost root of the system as a function of $a$ and $b$ for a fixed delay value of $\tau=0.2296$.}
\label{fig:surf}
\end{figure}

The optimization just presented assumed a fixed set of control gains and changed only the delay. By varying the three parameters it may be possible to move the rightmost root further to the left to increase the stability margin of the system. Figure \ref{fig:surf}, for example, presents the real part of the rightmost root of the system as a function of the control parameters $a$ and $b$ for a fixed delay value of $\tau=0.2296$. Here we observe that the control gain $a$ has a stronger effect than $b$ on the location of the rightmost root, and that the largest stability margin for this delay is obtained with control gains slightly larger than those selected. 

\begin{figure}[tb]
\centering
\includegraphics[scale=0.7]{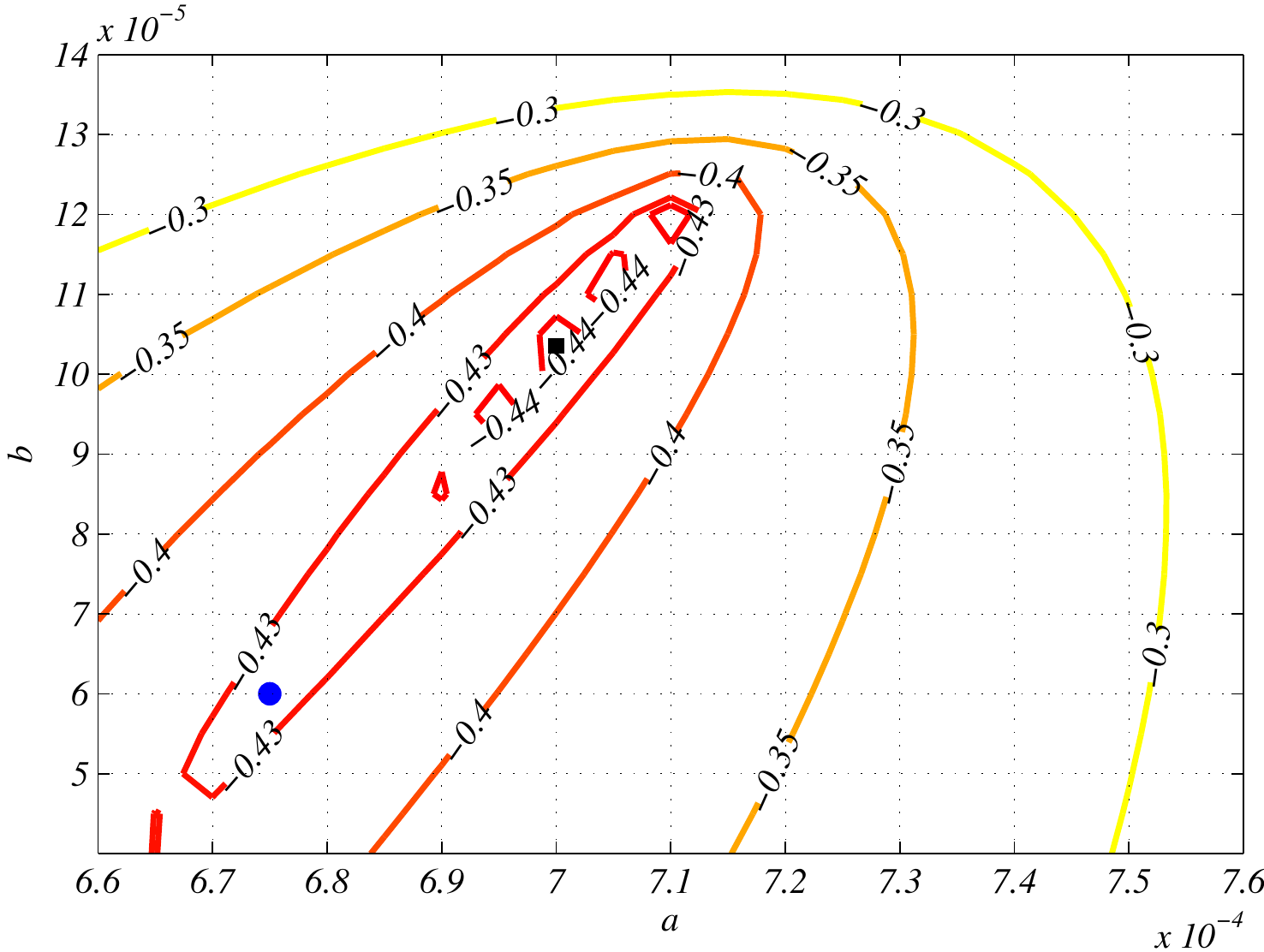}
\caption{Contour plot of the surface in figure \ref{fig:surf}}
\label{fig:contour}
\end{figure}

A better estimation of the optimal control gains for this delay value is made possible by observing Fig.~\ref{fig:contour}, the contour plot of the surface in Fig.~\ref{fig:surf}. For example, the control gains $a= 7\times 10^{-4}$ and $b= 1.03 \times 10^{-4}$, marked by a black square, locate the rightmost root of the system at $s= -0.4475+j1.7690$, giving a settling time a bit smaller than the one obtained with the current selection of control gains, marked with a blue dot in Fig.~\ref{fig:contour}. 

The knowledge of the stability regions in the domain of the delays allows the user to search for the best delay given a certain control gain selection. Since we are dealing with retarded systems, the stability region does not change dramatically as the parameters change and therefore a gradient descent algorithm can be implemented. With such systematic optimization routine, the best possible parametric $(a,b,\tau)$ combination for a given set of $\left(\sigma,\nu_1\right)$ parameters can be found using this approach.

\section{Closing Remarks}\label{sec:conc}
In this document we studied the problem of stabilization of pitch-flap motion of helicopter rotors using delayed feedback. We revisited a previously proposed control law and by using newer tools for the analysis of linear time-invariant, time delay systems we showed how the delays which make the system stable can be found. With his knowledge it is possible to find the largest possible stability margin of the system and we suggested a methodology to do so.

The Maple and Matlab codes used to create the plots in this document can be downloaded from \texttt{https://db.tt/MHvoLtdM}.

This document considered only the case of hovering flight, for which the system matrices are constant. In the case of forward flight, the dynamic model of the system becomes time periodic, and thus different tools need to be used for its analysis. This topic is a matter of further research. 

\bibliographystyle{ifacconf-hardvard}

\end{document}